\documentclass{article}

\usepackage{amsfonts,amsmath}
\usepackage{amssymb,latexsym}
\usepackage[dvips]{epsfig}
\usepackage{amscd}  

\title{\bf NilCoxeter algebras categorify the Weyl algebra}
\author{Mikhail Khovanov} 
\date{\today}

\newtheorem{prop}{Proposition}

\newtheorem{lemma}{Lemma}
\newtheorem{corollary}{Corollary}

\newcommand{\oplusop}[1]{{\mathop{\oplus}\limits_{#1}}}
\newcommand{\oplusoop}[2]{{\mathop{\oplus}\limits_{#1}^{#2}}}

\begin{document}
\maketitle
\baselineskip 12pt

\def\R{\mathbb R}
\def\N{\mathbb N}
\def\Z{\mathbb Z}
\def\Q{\mathbb Q}
\def\F{\mathbb F}
\def\S{\mathbb S}
\def\l{\lbrace}
\def\r{\rbrace}
\def\o{\otimes}
\def\D{\Delta}
\def\lra{\longrightarrow}

\def\C{\mc{C}}
\def\tr{\mathrm{tr}}

\def\a{\mathbf a} 
\def\l{\mathbf l} 
\def\H{\mathbb H} 
\def\M{\mc{M}}
\def\mC{\mc{C}}
\def\A{\mathcal{A}}
\def\Hom{\mbox{Hom}}
\def\wt{\widetilde}

\def\mc{\mathcal} 
\def\mf{\mathfrak}
\def\cC{{\mathcal{C}}}
\def\cR{{\mathcal{R}}}

\def\yesnocases#1#2#3#4{\left\{
\begin{array}{ll} #1 & #2 \\ #3 & #4 
\end{array} \right. }

\def\drawing#1{\begin{center} \epsfig{file=#1} \end{center}}

\def\hsm{\hspace{0.05in}}

\section{Introduction} 
\label{intro} 

In this paper we present an example of elaborate categorical structures 
hidden in very simple algebraic objects. We look at the algebra of 
polynomial differential operators in one variable $x,$ also known as  
the Weyl algebra, and its irreducible representation in the ring of 
polynomials $\Q[x].$ We construct an abelian category $\C$  whose 
Grothendieck group can be naturally identified with the ring of 
polynomials and define exact functors $F_X: \C \to \C$ and 
$F_D: \C \to \C$ such that 

(a) on the Grothendieck group $K(\C)$ of the category $\C$ functors 
$F_X$ and $F_D$ act as the multiplication by $x$ and differentiation, 
respectively, 

(b) there is a functor isomorphism $F_D F_X\cong F_X F_D\oplus 
\mbox{Id},$ which lifts the defining relation $\partial x = x \partial +1$ 
of the Weyl algebra, 

(c) functors $F_X$ and $F_D$ have nice adjointness properties: 
$F_X$ is left adjoint to $F_D$ and right adjoint to $F_D,$ twisted by 
an automorphism of $\C.$ 

The category $\C$ is the direct sum of categories $\C_n$ over 
all $n\ge 0,$ where $\C_n$ is the category of finite dimensional 
representations of the nilCoxeter algebra $A_n,$ which is 
generated by $Y_i, 1\le i \le n-1,$ subject to  relations 
$Y_i^2=0, Y_i Y_j=Y_j Y_i$ for $|i-j|>1$ and  
$Y_i Y_{i+1} Y_i = Y_{i+1} Y_i Y_{i+1}.$ 
The nilCoxeter algebra is the algebra of divided difference operators 
(see Macdonald [M], Fomin and Stanley [FS]) 
\begin{equation} 
Y_i f= \frac{f- s_i f}{x_{i+1}-x_i}, 
\end{equation}
where $f$ is a polynomial in $x_1, \dots , x_n$ 
and $s_i f(x_1, \dots , x_n) = f(x_1, \dots, x_{i+1}, x_i , \dots, 
x_n).$
Functors $F_X$ and $F_D$ are induction and restriction functors 
associated to the inclusion of algebras $A_n \hookrightarrow A_{n+1}.$ 
The following holds 

(d) functors $F_X$ and $F_D,$ restricted to each $C_n,$ are indecomposable. 

It seems likely that ours is the only example (up to obvious modifications 
of base change, etc.) of an abelian category 
$\C= \oplusop{n\ge 0} \C_n$ and 
exact functors $F_X$ and $F_D$ satisfying conditions (a)-(d) above.

In Section~\ref{Weyl} we study these and other properties of the category 
$\C$ and functors $F_X,F_D.$ In Section~\ref{bialgebra} we equip 
$\C$ with a bialgebra-category structure. Specifically, 
inclusions of algebras $A_n \otimes A_m \hookrightarrow A_{n+m}$ give rise 
to  
induction and restriction functors, which, when summed over all $n,m\ge 0,$ 
become functors $M: \C^{\o 2} \lra \C$ and $\Delta: \C \to \C^{\o 2}$ 
between $\C$ and its second tensor power $\C^{\otimes 2}.$ 
These functors are exact, boast neat adjointness properties 
and on the Grothendieck group descend to the 
multiplication and comultiplication in the commutative, cocommutative 
Hopf algebra $\Q[x]$ of polynomials in one variable.  
The associativity relation for the  multiplication, coassociativity of the  
comultiplication and the consistency relation between the multiplication and 
comultiplication become isomorphisms of functors. We check that these 
isomorphisms satisfy the coherence relations of Crane and Frenkel [CF] 
for a bialgebra-category.  

In Section~\ref{graded} we sketch how working with graded modules and 
bimodules over the nilCoxeter algebra yields a categorification of 
the quantum Weyl algebra and of 
the quantum deformation of the Hopf algebra $H.$  
The grading shift automorphism in the category of graded modules
descends to a map of Grothendieck groups which we interpret as the 
multiplication by $q.$ 

In Section~\ref{wreath} we present a simple generalization of our 
construction to cross-products and provide a short 
comment on the relation of our work to Ariki's 
realization [A] of highest weight modules over affine Lie algebras. 

\vspace{0.1in} 

{\it Acknowledgements: } This work was done during my very nice stay 
at the Institute for Advanced Study. I would like to thank the Institute 
and the NSF for supporting me with grants DMS 97-29992 and DMS 96-27351. 

\section{The Weyl algebra and bimodules over the nilCoxeter algebra} 
\label{Weyl} 

\subsection{The Weyl algebra} 
\label{the-Weyl} 

The Weyl algebra $W$ is the algebra of differential operators with polynomial 
coefficients in one variable. For our purposes define $W$ as the 
algebra over $\Z$ with the unit $1,$ generators $x, \partial$ and defining 
relation $\partial x= x\partial +1.$ Let $R_{\Q}$ be the $\Q$-vector 
space $R_{\Q}$ spanned by $x^0,x^1, x^2,\dots$ .
$W$ acts on $R_{\Q}$ via $x\cdot x^i= x^{i+1}$ and 
$\partial \cdot x^i = i x^{i-1}.$ Abelian subgroups $R$ and $R'$ of $R_{\Q},$ 
generated by $\{ x^i/i! \}_{i=0}^{\infty}$ and 
$\{ x^i\}_{i=0}^{\infty},$ respectively, 
are $W$-submodules of $R_{\Q},$ i.e. the action 
of $x$ and $\partial$ has integral coefficients in each of these two bases. 

The Weyl algebra has an antiinvolution $\tau: W\to W$ with 
\begin{equation} 
\tau(x)= \partial, \hspace{0.15in} \tau(\partial) =x
\hspace{0.15in} \mbox{ and } \hspace{0.15in} 
\tau(ab) = \tau(b)\tau(a) \hspace{0.1in} 
\mbox{ for } \hspace{0.1in} a,b\in W.
\end{equation} 
Let $(,)$ be the symmetric bilinear form on $R_{\Q}$ defined by 
$(x^i,x^j) = \delta_{i,j} i!.$ This form is $\tau$-invariant: 
\begin{equation} 
( y a,b) = (a,\tau(y)  b) \hspace{0.1in} \mbox{for} \hspace{0.1in} 
 y \in W, \hspace{0.1in} a,b \in R_{\Q}, 
\end{equation} 
and it restricts to an integer valued bilinear product 
$(,): R' \times R \to \Z.$ 

\subsection{The nilCoxeter algebra} 
\label{nilCoxeter} 

Let $A_n$ be the unital algebra over $\Q$ generated by 
$Y_1, \dots , Y_{n-1}$ with defining relations 
\begin{equation} 
\label{three} 
\begin{array}{rcl} 
  Y_i^2 & = & 0  \\
        &   &     \\
  Y_i Y_j & = & Y_j Y_i \hspace{0.2in} |i-j|>1 \\
        &   &       \\
  Y_i Y_{i+1} Y_i & = & Y_{i+1} Y_i Y_{i+1}. 
\end{array} 
\end{equation} 
Fomin and Stanley [FS] call $A_n$ {\it the nilCoxeter algebra}. 
It originally appeared in the work of Bernstein, Gelfand and Gelfand [BGG] on 
the cohomology of flag varieties and was later investigated and 
generalized in various ways by Lascoux and 
Sch\"{u}tzenberger [LS], Macdonald [Mc], Kostant and Kumar [KK],
Fomin and Stanley [FS] and others.  
Note that if we change the first relation in (\ref{three}) to $Y_i^2=1,$ 
we obtain the group algebra of the symmetric group, which is, indeed, 
closely related to the nilCoxeter algebra: 

\begin{prop} \label{nice-prop} The algebra $A_n$ is isomorphic 
to the algebra which is spanned over $\Q$ by $Y_{w},$ as $w$ ranges over 
 elements of the symmetric group $\S_n,$ with the multiplication
\begin{equation} 
\begin{array}{rcl}  
Y_{w_1} Y_{w_2} & = & Y_{w_1 w_2}\hspace{0.1in} \mbox{ if } 
 \hspace{0.1in} l(w_1w_2)= l(w_1)+ l(w_2),  \\
                &   &       \\
Y_{w_1} Y_{w_2} & = &  0 \hspace{0.35in} \mbox{ otherwise,}
\end{array}
\end{equation}  
where $l(w)$ is the standard length function on the symmetric group, the 
number of inversions created by $w.$  
The isomorphism is given by sending the generator $Y_i$ of $A_n$ 
to $Y_{s_i}$ where $s_i= (i,i+1)$ is the transposition of $i$ and $i+1.$ 
\end{prop} 

In particular, $A_n$ has dimension $n!.$ Note that $A_0= A_1=\Q.$ 
Introduce a trace map $\tr_n: A_n\to \Q$ by 
\begin{equation} 
\begin{array}{rcl}
\tr_n(Y_{w_0}) & = & 1 \hspace{0.1in} \mbox{ where } w_0
 \mbox{ is the longest permutation,} \hspace{0.1in} w_0(i)= n-i,   \\
            &   &    \\
\tr_n(Y_w)  & = & 0 \hspace{0.1in} \mbox{ if } w\not= w_0.
\end{array}  
\end{equation}

\begin{prop} The trace map $\tr_n$ is nondegenerate and makes $A_n$ into 
a Frobenius algebra.  
\end{prop} 

{\emph Proof:} When we say that the trace map is nongenerate we mean that for 
each $y\in A_n, y\not=0$ we can find $y'\in A_n$ such that $\tr_n(y y')=1.$ 
Algebras with a nondegenerate trace map are called Frobenius algebras. 
The basic properties of the length function in the symmetric group 
imply that $\tr_n$ is nondegenerate. $\square$ 

For more information about Frobenius algebras we refer the reader 
to Yamagata [Y] and references therein. 

Let $B_1,B_2$ be finite-dimensional algebras over a field $k$ and 
$N$ a finite-dimensional  $(B_1,B_2)$-bimodule. Then $N^{\ast}= \Hom_k(N,k)$ 
is naturally a $(B_2, B_1)$-bimodule. The duality functor $\ast$ is a 
contravariant equivalence between categories of finite-dimensional 
$(B_1, B_2)$-bimodules  and $(B_2,B_1)$-bimodules. When $B_2=k,$ the duality 
functor $\ast$  is a contravariant 
 equivalence between categories of left and right 
$B_1$-modules. 

If $B_1$ has an automorphism $\nu,$ we can use it to twist the 
left action of $B_1$ on a bimodule $N:$ for $y\in B_1$ and $t\in N$ 
the twisted left action of $B_1$ is $\nu(y)t.$ We will denote the resulting 
bimodule by $_{\nu}N.$ An automorphism $\nu$ of $B_2$ allows to twist the 
right action of $B_2$ on $N,$ we denote the resulting bimodule by 
$N_{\nu}.$ 

Any algebra $B$ is a bimodule over itself in the obvious way. 
Denote by $\psi_n$ the involution of $A_n$ which takes 
$Y_i$ to  $Y_{n-i}$ and by $A_n^{\psi}$ the algebra $A_n$ as a 
bimodule over itself with the right action twisted by $\psi_n.$  
Let $1_n^{\psi}\in A_n^{\psi}$ 
be the image of $1\in A_n$ under the isomorphism 
$A_n \stackrel{\cong}{\longrightarrow} A_n^{\psi}$ of right $A_n$-modules, 
so that $1_n^{\psi} Y_i = Y_{n-i} 1_n^{\psi}.$

\begin{prop} 
\label{iso-iso}
$A_n$-bimodules $A_n^{\ast}$ and $A_n^{\psi}$ are isomorphic. 
\end{prop}  

{\emph Proof: } For $w\in \S_n$ let $Y_w^{\ast}\in A_n^{\ast}$ be the 
functional 
\begin{equation} 
Y_w^{\ast}(Y_{\sigma}) = 
\begin{cases}   1 & \text{if $\sigma w= w_0$}, \\
                0 & \text{otherwise} 
\end{cases} 
\end{equation}  
It is easy to check that the map $A_n^{\ast}\to A_n^{\psi}$ given by 
$Y_w^{\ast} \longmapsto Y_w 1_n^{\psi}$ is a 
bimodule isomorphism (use that $s_i w_0 = w_0 s_{n-i}$). $\square$ 

\begin{corollary} 
The algebra $A_n$ is injective as a left and right $A_n$-module.  
\end{corollary} 

This follows from either of the last two propositions. $\square$ 

\vspace{0.1in} 

Note that the trace $\mbox{tr}_n$ is quasi-symmetric w.r.t. $\psi_n$: 
\begin{equation} 
\mbox{tr}_n ( a b) = \mbox{tr}_n(\psi_n(b) a)\hspace{0.1in} 
\mbox{ for } \hspace{0.1in} a,b\in A_n.
\end{equation} 
In the terminology of Frobenius algebras, $\psi_n$ is the Nakayama 
automorphism associated with $\mbox{tr}_n$ (see [Y], Section 2.1).

\subsection{Bimodules} 
\label{bimod} 

Denote by  $\chi_n$ the algebra map $A_n \to A_{n+1}$ which sends each 
$Y_i$ to $Y_i.$ Proposition \ref{nice-prop} implies that $\chi_n$ is 
injective. 
The inclusion $\chi_n: A_n \to A_{n+1}$ induces a left and right 
$A_n$-module structure on $A_{n+1}.$ The left $A_n$ module structure 
on $A_{n+1}$ commutes with the right $A_{n+1}$-module structure on 
$A_{n+1},$ the latter coming from the right action of $A_{n+1}$ on itself. 
Thus, $A_{n+1}$ is an $(A_n,A_{n+1})$-bimodule in a natural way, and 
we denote this bimodule by $D_{n+1}.$ Similary, we get an  
$(A_{n+1},A_n)$-bimodule structure on $A_{n+1}$ and denote this bimodule 
by $X_n.$ 

\begin{prop} Bimodules $X_n$ and $D_n$ are left and right projective. 
\end{prop} 
{\emph Proof:} Bimodule $X_n$ is free of rank $1$ as a left $A_{n+1}$-module 
and, thus, left projective. As a right $A_n$-module, it is free of 
rank $n+1$ and spanned by $1, Y_{n}, Y_{n-1}Y_n, Y_{n-2}Y_{n-1}Y_n, 
\dots, Y_1Y_2\dots Y_n$
(since any element $s$ of the symmetric group $\S_{n+1}$ admits 
a unique decomposition $s=  s_is_{i+1}\dots s_n s'$ with $s'\in \S_n$ 
 and some $i, 1\le i \le n+1$). Hence, $X_n$ is projective as a 
right $A_n$-module. The same argument works for $D_n.$ 
$\square$ 

\begin{prop} 
\label{n-isom-prop} 
For each $n,$ there is an isomorphism of $A_n$-bimodules 
\begin{equation} 
\label{bimod-eqn} 
     D_{n+1} \o_{A_{n+1}} X_n \cong  A_n \oplus  (X_{n-1} \o_{A_{n-1}} D_n) , 
\end{equation} 
where $A_n$ is equipped with the standard bimodule structure. 
\end{prop} 

{\emph Proof: } The left hand side of (\ref{bimod-eqn}) is isomorphic 
to $A_{n+1},$ considered as an $A_n$-bimodule via $\chi_n.$ The right 
hand side is the direct sum of $A_n$ and  $A_n \o_{A_{n-1}} A_n.$ 
We have maps $m_1, m_2$ of $A_n$-bimodules 
\begin{equation} 
m_1: A_n \to A_{n+1} \hspace{0.1in} \mbox{ and } \hspace{0.1in} 
m_2: A_n \o_{A_{n-1}} A_n \to A_{n+1}
\end{equation} 
which are uniquely determined by $m_1(1) = 1$ and $m_2(1\o 1) = Y_n.$ 
For $w\in \S_{n+1},$ the element $Y_w$ of  $A_{n+1}$ lies in 
$m_1(A_n)$ iff $w(n+1)= n+1.$ If $w(n+1)\not= n+1,$  we can write 
$w= y s_n z$ for $y,z\in \S_n,$ so that $Y_w= m_2(Y_y\o Y_z).$ 
Therefore, $m_1$ and $m_2$ define an $A_n$-bimodule isomorphism 
\begin{equation} 
A_n \oplus (A_n \o_{A_{n-1}} A_n) \stackrel{\cong}{\longrightarrow} A_{n+1}.  
\end{equation} 
$\square$ 

\vspace{0.2in} 

Let $A= \oplusoop{n=0}{\infty} A_n$ be the direct sum of algebras $A_n$ 
over all $n.$ Algebra $A$ does not have a unit, instead it has an 
infinite system of pairwise orthogonal idempotents $1\in A_n, n\ge 0.$  

An $(A_n, A_k)$-bimodule $N$ is naturally a bimodule over the algebra $A.$ 
 Namely, for $x\in A_i, i\not= n$ we set $xN=0$ and let $Nx=0$ if 
$x\in A_i , i\not= k.$ 
In this way, bimodules $X_n$ and $D_n,$ as we sum over all $n,$ 
give rise to $A$-bimodules $X= \oplusop{n\ge 0} X_n$ and 
$D= \oplusop{n\ge 0} D_n.$ 
We can rephrase Proposition \ref{n-isom-prop} as 

\begin{prop} \label{nat-prop} 
There is a natural isomorphism of $A$-bimodules 
\begin{equation} 
 D \o_{A} X \cong  A \oplus  (X \o_{A} D). 
\end{equation} 
\end{prop} 

This is, of course, a bimodule version of the Weyl algebra relation 
$\partial x = x\partial + 1.$ The generators $x$ and $\partial$ 
of the Weyl algebra become the  
bimodules $X$ and $D,$ the  product in the Weyl algebra becomes the tensor 
product of bimodules,  addition becomes the direct sum and $1$ becomes the 
identity bimodule $A.$ 

In the rest of this section we continue in the similar fashion, 
interpreting other structures of the Weyl algebra and its polynomial 
representation in the framework of nilCoxeter algebras. 

\subsection{Categories and functors} 
\label{cat-and-fun} 

Let $\C_n$ be the category of finite-dimensional unital left 
$A_n$-modules, and  let $\C= \oplusoop{n=0}{\infty} \C_n.$ 
The category $\C$ can be viewed as 
the full subcategory of the category of finite-dimensional left $A$-modules, 
which consists of $A$-modules $N$ with $AN = N$ and $A_nN=0$ for 
large enough $n.$ 

An $A$-bimodule $T$ is called \emph {small} if it preserves 
the category $\C,$ 
i.e., for any $N\in \C,$ the module $T \o_{A} N$ is in $\C.$  
Denote by $F_T$ the functor of tensoring with $T.$ 
We can reformulate Proposition \ref{nat-prop} as saying that there is 
a canonical isomorphism of functors 
\begin{equation} 
\label{the-equation}
F_D F_X \cong F_X F_D\oplus \mbox{Id}_{\C}
\end{equation} 

The Grothendieck group $K(\mc{U})$ of an abelian category $\mc{U}$ 
is the group generated by symbols $[N]$ for all objects $N$ of $\mc{U}$ 
subject to relations $[N_2]= [N_1]+ [N_3]$ whenever there is a 
short exact sequence $0 \to N_1\to N_2\to N_3\to 0.$ 
The Grothendieck group of $\C_n$ is isomorphic to $\Z$ and generated 
by $[L_n],$ where $L_n$ is the simple $A_n$-module ($L_n$ is 
uniquely defined, up to an isomorphism). We will identify $K(\C_n)$ 
with the abelian subgroup of $R$ generated by $x^n/n!,$   by sending 
$[L_n]$ to $x^n/n!.$ 

Since $K(\C)= \oplusop{n\ge 0} K(\C_n),$ 
the Grothendieck group of $\C$ is canonically identified with 
the abelian group $R,$ so that from now on we will consider $[N],$ for 
any object $N$ of $\C,$ as an element of $R.$ 
The indecomposable projective module 
in $\C_n$ (which we denote $P_n$) is mapped to $x^n$:  
\begin{equation} 
[P_n] = x^n, \hspace{0.15in} [L_n]= \frac{x^n}{n!}. 
\end{equation}  
We interpret the bilinear form $(,): R'\times R\to \Z$ via 
the $\Hom$ bifunctor:  if $P$ is a projective object of $\C$  and $N$ 
is any object, then 
\begin{equation} 
\label{its-bilinear} 
\dim_{\Q}(\Hom_{\C}(P,N)) = ([P], [N]) 
\end{equation} 
Note that we need $P$ to be projective, otherwise the dimension function 
on the left hand side will not be additive in $[N].$ Observe also that 
the form $(,)$ takes values in $\Z$ when restricted to $R'\times R,$ but 
is fractional when restricted to $R\times R.$

\hspace{0.15in}

Bimodules $X$ and $D$ are right projective, so that the functors of tensoring 
with them are exact and induce maps $x$ and $\partial$ of the Grothendieck 
group $R= K(\C)$: 
\begin{equation} 
[X\o_{A} N] = x [N], \hspace{0.15in} 
[D\o_{A} N]=  \partial [N] \hspace{0.15in} \mbox{ for } 
\hspace{0.15in} N\in \mbox{Ob }\C. 
\end{equation}

The functor of tensoring with $X_n$ is the induction functor from 
$A_n$-modules to $A_{n+1}$-modules, while tensoring with $D_{n+1}$ 
is the restriction functor from $A_{n+1}$-modules to $A_n$-modules. 
Since the induction is left adjoint to the restriction, we conclude that 
$F_{X_n}$ is left adjoint to $F_{D_{n+1}}$ and $F_X$ 
{\it is left adjoint to}  
$F_D,$ i.e., there is a bifunctor isomorphism 
\begin{equation} 
\Hom_{\C}(F_X?, ?) \cong \Hom_{\C}(?, F_D ?). 
\end{equation} 
This isomorphism can be interpreted as the lift of the equality 
$(x a,b) = (a, \partial b)$ for $a,b\in R_{\Q},$ since we just established 
that the $\Hom$-bifunctor lifts the bilinear form $(,)$  (formula 
(\ref{its-bilinear})). 

Note that $(\partial a,b) = (a, xb),$ so that a natural 
guess says that $F_X$ is not only left but 
also right adjoint to $F_D.$ This is false,  
but not far from the truth. Consider the bimodule $A_n^{\psi},$ which 
was defined in Section \ref{nilCoxeter}. Denote by $A^{\psi}$ the 
$A$-bimodule which is the direct sum of $A_n^{\psi}$ over all $n,$ and 
by $\Psi: \C \to \C$ the functor $F_{A^{\psi}}$ 
of tensoring with $A^{\psi}.$ Since $\psi_n$ is an involution, 
$\Psi^2 \cong \mbox{Id}_{\C}.$ 

\begin{prop} The functor $F_X$ is right adjoint to 
$\Psi F_D \Psi.$ 
\end{prop} 

\emph{Proof: } We use the following 

\begin{lemma}
\label{twist} Let $B_1,B_2$ be Frobenius algebras over a field $k$ and 
$\nu_1,\nu_2$ be Nakayama automorphisms of $B_1,B_2.$ Suppose that 
$N$ is a finite-dimensional $(B_1,B_2)$-bimodule which is projective 
as a left $B_1$-module and as a right $B_2$-module. Then the functor 
\begin{equation} 
 N \o_{B_2} ? : B_2\mathrm{-mod} \to B_1\mathrm{-mod}
\end{equation} 
of tensoring with $N$ has the right adjoint functor 
\begin{equation} 
 (N^{\ast})_{\nu_1}\o_{B_1} ? : B_1\mathrm{-mod} \to B_2\mathrm{-mod} 
\end{equation} 
(here $(N^{\ast})_{\nu_1}$ is the dual of $N,$ with the right $B_1$-action 
twisted by $\nu_1$) and the left adjoint functor 
\begin{equation} 
_{\nu_2^{-1}}(N^{\ast})\o_{B_1} ? : B_1\mathrm{-mod} \to B_2\mathrm{-mod} 
\end{equation} 
\end{lemma} 

In the case when 
 $B_1$ and $B_2$ are symmetric algebras (i.e. $\nu_1,\nu_2$ are identity 
maps), 
rather than just Frobenius algebras, 
this lemma is proved in Rickard [R], Corollary 9.2.4. The same proof works 
for Frobenius algebras. $\square$ 

Applying this lemma to the $(A_{n+1}, A_n)$-bimodule $X_n$ proves 
the proposition (the Nakayama automorphism of $A_n$ is $\psi_n,$ hence 
the conjugation by $\Psi$ in the second adjointness isomorphism). 
$\square$

\vspace{0.15in} 

The algebra 
$W$ has a $\Q$-vector space basis $\{ x^m \partial^n\}$ for 
$n,m\ge 0.$ We will call this basis 
 the canonical basis of $W.$ A product of two elements of the canonical 
basis decomposes as a linear combination of canonical basis vectors with 
nonnegative integral coefficients. 
This basis can be interpreted in terms of indecomposable bimodules. Namely,  
the $(A,A)$-bimodule $X^{\o m}\o D^{\o n}$ (all tensor products are 
over $A$), 
which in our theory is naturally associated to $x^m \partial^n,$   
is the direct sum of  $(A_{m+k-n}, A_k)$-bimodules 
$A_{m+k-n} \o_{A_{k-n}} A_k,$ over all $k\ge n,$  and we have  

\begin{prop} The $(A_{m+k-n},A_k)$-bimodule $A_{m+k-n} \o_{A_{k-n}} A_k$ 
is indecomposable. 
\end{prop} 

\emph{Proof:} An exercise. $\square$

\subsection{Contravariant duality} 
\label{dualities}

Denote by $\C_n^r$ the category of finite-dimensional right $A_n$-modules 
and by $\C^r$ the direct sum of categories $\C_n^r$ over all $n\ge 0.$ 
The duality functor $N^{\ast} = \mbox{Hom}(N,k),$ defined in 
Section~\ref{nilCoxeter} for bimodules, will be considered in 
this section as a contravariant functor from $\C$ to $\C^r.$ 

Let  $u$ be the antiinvolution of $A_n$ which takes $Y_i$ to $Y_i.$ 
It induces an equivalence of categories of left and right $A_n$-modules. 
As we sum over all $n,$ we obtain an  equivalence of
categories $U: \C^r \to \C.$ 

The functor $\Omega = U \ast$ is a contravariant equivalence $\C \to \C.$ 

\begin{prop} There are natural isomorphisms of functors
\begin{eqnarray} 
  \label{square-is-one} \Omega^2      & \cong  & \mathrm{Id}_{\C}   \\
  \label{box-psi} \Omega \Psi    & \cong  & \Psi \Omega  \\
  \label{box-X}   \Omega \Psi F_X     & \cong  & F_X \Omega \Psi \\
  \label{box-D}   \Omega F_D     & \cong  & F_D \Omega 
\end{eqnarray} 
\end{prop} 

\emph{Proof}  $\Omega^2\cong \mbox{Id}_{\C},$
since $u$ is an antiinvolution, and $\Omega \Psi \cong \Psi \Omega,$ since
$u\psi_n = \psi_n u.$ Isomorphism (\ref{box-X}) is a corollary of  

\begin{lemma} 
\label{dual} 
\begin{enumerate} 
\item There is an isomorphism of bimodules 
\begin{equation} 
X_n^{\ast} \cong  D_{n+1} \o_{A_{n+1}} A_{n+1}^{\psi}
\label{x-and-d} 
\end{equation} 
\item There are  isomorphisms, functorial in $N\in \C_n,$
\begin{eqnarray}
 (A_n^{\psi}\o_{A_n} N)^{\ast} & \cong & N^{\ast}\o_{A_n} A_n^{\psi} \\ 
 (X_n \o_{A_n} N )^{\ast} & \cong & N^{\ast}\o_{A_n} \widetilde{D}_{n+1}, 
\label{x-and-m}
\end{eqnarray} 
where $\widetilde{D}_{n+1}$ is the $(A_n,A_{n+1})$-bimodule 
$A_n^{\psi}\o_{A_n} D_{n+1}\o_{A_{n+1}} A_{n+1}^{\psi}.$ 
\item There are functorial in $N\in C_n^r$ isomorphisms 
\begin{equation} 
U(N\o_{A_n} D_{n+1}) \cong X_n \o_{A_n} U(N), 
\hspace{0.2in} 
U(N\o_{A_n} A_n^{\psi}) \cong A_n^{\psi}\o_{A_n} U(N). 
\end{equation} 
\end{enumerate} 
\end{lemma} 

Statement 1 of the lemma follows from  Proposition \ref{iso-iso}.
Let us now prove (\ref{x-and-m}). If  $N=P_n,$
the indecomposable projective $A_n$-module, the  
isomorphism (\ref{x-and-m}) follows from (\ref{x-and-d}). 
Moreover, (\ref{x-and-d}) also implies that 
(\ref{x-and-m}) is functorial relative to $A_n$-module maps $P_n \to P_n$
(here $N= P_n$). For an arbitrary $N,$ represent $N$ as the cokernel of 
a map of projective modules:  $P_n^{\oplus a}\to P_n^{\oplus b} 
\to N \to 0.$ Applying the 
functors on the left and right hand sides of (\ref{x-and-m}) to each 
term of this 
exact sequence, and using the exactness of tensoring with $X_n$ and 
$\widetilde{D}_{n+1},$ we conclude that (\ref{x-and-m}) holds for 
any $N,$ functorially in $N.$ Other statements of the lemma can be proved 
in a similar or  easier fashion. $\square$ 

Armed with Lemma~\ref{dual}, we compute, for $N\in \C,$ 
\begin{equation} 
\Omega \Psi F_X (N)  = \Psi \Omega ( X \o N) = \Psi U( N^{\ast} \o 
\widetilde{D})  = \Psi^2 X\o (\Psi U(N^{\ast})) = F_X \Psi \Omega (N). 
\end{equation} 

Isomorphism (\ref{box-D}) is adjoint to (\ref{box-X}). $\square$


\subsection{The integral} 
\label{integral} 

We can next ask about the meaning of  
the indefinite integral in our model. The formula 
\begin{equation} 
\label{integrals} 
\int x^n = \frac{x^{n+1}}{n+1}
\end{equation} 
 suggests to look for an 
exact functor from $\C_n$ to $\C_{n+1}$ which takes the projective module 
$P_n$ to a module which is $n+1$ times ``smaller'' than the projective 
module $P_{n+1}$ (since in our correspondence the image of the projective 
module $P_i$ in the ring of polynomials is $x^i$). Let $I_n$ be the 
$(A_{n+1},A_n)$-bimodule, which is isomorphic to $A_n$ as the right 
$A_n$-module, and the left $A_{n+1}$-action is via the homomorphism 
of algebras $t_{n+1}: A_{n+1}\to A_n, t_{n+1}(Y_i)= Y_i$ for $i<n$ and 
$t_{n+1}(Y_n)=0.$ Since $I_n$ is projective as a right $A_n$-module, 
the functor $F_{I_n}: \C_n \to \C_{n+1}$ of tensoring an $A_n$-module 
with $I_n$ is exact. Moreover, $F_{I_n}$ takes the indecomposable
projective module 
$P_n$ to a module of dimension $n!,$ while the projective generator 
$P_{n+1}$ of $\C_{n+1}$ has dimension $(n+1)!,$ so that the desired 
relation holds: $[ I_n \o_{A_n} N ] = \int  [N]$ for   $N\in \mbox{Ob}(\C_n).$ 
To formulate this relation without the index $n,$ we form 
$I= \oplusoop{n=0}{\infty} I_n,$ the $A$-bimodule which is the 
direct sum of $I_n$ over all $n.$ 
Then we have
\begin{equation} 
[I\o_{A}N] = \int  [N]\hspace{0.15in} \mbox{for all} 
\hspace{0.15in} N \in \mbox{Ob} \C. 
\end{equation}  

The following result is obvious: 
\begin{prop} There are  bimodule isomorphisms 
$D_{n+1}\o_{A_{n+1}} I_n \cong  A_n$ and 
$D \o_{A} I \cong A.$ 
\end{prop}   
This isomorphism can be considered as a categorification of the 
formula $d \int f(x) = f(x),$  
for a polynomial 
function $f(x).$ On the other hand, we don't get to categorify 
the formula $\int d f(x) = f(x),$ for  
 $I\o_{A} D$ is not isomorphic to $A$ as an $A$-bimodule.

\subsection{Multiplication and the Leibniz rule} 
\label{Leibniz} 

Let $\gamma_{n,m}$ be the algebra homomorphism $A_n \o A_m\to A_{n+m}$ 
given by 
\begin{equation} 
\gamma_{n,m}(Y_i\o 1) = Y_i, \hspace{0.2in} 
\gamma_{n,m} (1\o Y_i) = Y_{n+i}. 
\end{equation} 
$\gamma_{n,m}$ is injective and induces a bifunctor, denoted 
$J_{n,m},$ from the product $\C_n \times \C_m$ of categories $\C_n$ 
and $\C_m$ to $\C_{n+m}$: 
\begin{equation} 
J_{n,m}(N_1,N_2) =   A_{n+m}\o_{(A_n\o A_m)} (N_1\o N_2)
\hspace{0.15in} \mbox{ for } \hspace{0.15in} N_1\in \C_n, N_2\in \C_m.
\end{equation} 
Denote by $J$ the bifunctor $\C\times \C\to \C,$ which is the direct 
sum of $J_{n,m}$ over all $n,m\ge 0.$ 

\begin{prop} 
\begin{enumerate} 
\item Bifunctor $J$ is biexact. 
\item There is a functorial isomorphism 
$F_D \circ J(N_1,N_2) \cong J(N_1, F_D N_2)\oplus J(F_D N_1,N_2),$ satisfying 
the natural consistency relation for the decomposition of 
$F_D \circ J(J(N_1,N_2), N_3).$ 
\end{enumerate} 
\end{prop} 
We omit the proof. $\square$

Since $J$ is biexact, it induces a map of Grothendieck groups 
$K(\C)\times K(\C) \to K(\C),$ which is just the multiplication in the 
ring of polynomials.  Part 2 of the proposition 
is a functor version of the Leibniz rule $\partial(ab)= (\partial a)b+ 
a (\partial b).$

\section{The bialgebra-category structure of $\C$} 
\label{bialgebra} 

\subsection{Multiplication and comultiplication} 

The algebra of polynomials $R_{\Q}= \Q[x]$ has a comultiplication 
$c(x) = x\o 1+ 1\o x$ which makes $R_{\Q}$ into a bialgebra. The subring 
$R$ of $R_{\Q}$ is stable under the comultiplication and has a structure 
of a bialgebra over $\Z.$ We will  explains in detail how to lift 
the bialgebra structure from $R$ to the category $\C.$ 

Let ${\bf n}= (n_1, \dots, n_i)$ be an ordered 
$i$-tuple of nonnegative integers. 
Let $A_{\bf n}= A_{n_1}\o \dots\o A_{n_i}$ and denote by $\C_{\bf n}$ the 
category of finite dimensional left $A_{\bf n}$-modules. Let $\C^{\o i}$ 
the direct sum of categories $\C_{\bf n}$ over all $i$-tuples ${\bf n}.$ 

Algebra homomorphisms $\gamma_{n,m}: A_n \o A_m \to A_{n+m}, $ summed over 
all $n$ and $m,$ define induction and restriction functors: 
\begin{equation} 
M: \C^{\o 2} \to \C, \hspace{0.15in}\Delta: \C \to \C^{\o 2}. 
\end{equation} 

Note that the Grothendieck group of $\C^{\o i}$ is naturally isomorphic 
to the $i$-th tensor power of $K(\C).$ 
The symmetric group $\mathbb{S}_i$ acts on the set of $i$-tuples by 
permutations of terms. This action induces an action of $\mathbb{S}_i$ on 
the category $\C^{\o i}.$ We denote by $S_{j,j+1}$ the action of the 
transposition $(j, j+1).$ 

If a functor $G_k: \C^{\o i_k} \to \C^{\o j_k}, k=1,2$ is given by 
tensoring with a bimodule $W_k,$ denote by $G_1\o G_2$ the functor 
$\C^{\o (i_1+ i_2) } \to \C^{\o (j_1+ j_2)}$ of tensoring with 
the bimodule $W_1\o_{\Q} W_2.$ 

\begin{prop} 
\begin{enumerate} 
\item $M$ is left adjoint to $\Delta$ and  right adjoint to $S_{12} \Delta.$  
\item There are functor isomorphisms 
\begin{eqnarray} 
 M S_{12} & \cong & \Psi M \Psi^{\o 2}  \label{quasi-com-mult} \\
 S_{12} \Delta & \cong & \Psi^{\o 2} \Delta \Psi \label{quasi-com-delta}  
\end{eqnarray} 
\item 
Functors $M$ and $\Delta$ are exact and on the Grothendieck 
group  descend to the multiplication and comultiplication in the 
bialgebra $K(\C).$ 
\end{enumerate} 
\end{prop} 

\emph{Proof} Part 2 of the proposition follows from an obvious computation 
with bimodules. Next, 
$M$ is induction and $\Delta$ is restriction, thus, $M$ is 
left adjoint to $\Delta.$ $A_{n+m}$ is projective as left or 
right $A_n\o A_m$-module, so we can apply Lemma~\ref{twist}
and conclude that $M$ is right adjoint to $\Psi^{\o 2} \Delta \Psi.$ 
Together with the isomorphism (\ref{quasi-com-delta}), this implies that 
$M$ is right adjoint to $S_{12}\Delta.$  
Since $M$ and $\Delta$ each have left and right adjoints, they are 
exact. $\square$

\begin{prop} 
\label{funct-iso} 
\begin{enumerate} 
\item There are functor isomorphisms 
\begin{eqnarray} 
M (M \o \mathrm{Id}) &  \cong & M (\mathrm{Id} \o M), \label{iso-ass} \\
(\Delta \o \mathrm{Id}) \Delta & \cong & (\mathrm{Id}\o \Delta) \Delta,  
\label{iso-co}  \\
\label{consistency} 
\Delta M & \cong &  M^{\o 2} S_{23} \Delta^{\o 2}. 
\end{eqnarray} 
\end{enumerate} 
\end{prop}

\emph{Proof } Functors $M(M\o \mbox{Id})$ and $M(\mbox{Id}\o M),$ restricted 
to $\C_n\o \C_m \o \C_k,$ are canonically isomorphic to the functor of 
tensoring with $A_{n+m+k},$ considered as a left $A_{n+m+k}$-module and 
a right $A_n \o A_m \o A_k$-module. Hence the functor isomorphism 
(\ref{iso-ass}). The same argument works for (\ref{iso-co}). To prove 
(\ref{consistency}), note that both sides of it decompose as direct 
sums of functors $\C_n \o \C_m \to \C_k \o \C_l,$ over all quadruples 
$(n,m,k,l)$ such that $n+m = k+l.$ The left hand side of 
(\ref{consistency}), as a functor $\C_n \o \C_m\to \C_k \o \C_l,$ 
is naturally isomorphic to the functor of tensoring with $A_{n+m},$ 
considered as a left $A_k \o A_l$ and a right $A_n \o A_m$-module, via 
algebra homomorphisms $\gamma_{k,l}$ and $\gamma_{n,m}.$

\begin{lemma} Let $w_1, \dots, w_p,$ where $p = \mbox{min}(n,m,k,l),$ 
be minimal length representatives of double cosets 
$\S_k \times \S_l \setminus \S_{k+l} / \S_n \times \S_m.$ Then 
$A_{n+m}$ is isomorphic, as an $(A_k\otimes A_l, A_n \otimes A_m)$-bimodule, 
to the direct sum of subbimodules of $A_{n+m},$ spanned by $Y_{w_1}, 
\dots , Y_{w_p}.$ 
\end{lemma} 

We omit the proof of the lemma. $\square$ 

The right hand side of (\ref{consistency}), as a functor 
$\C_n \o \C_m \to \C_k \o \C_l,$ is isomorphic to the direct sum 
(over all admissible $r$) of the following functors: restrict from $A_n\o A_m$ 
to $A_r\o A_{n-r}\o A_{k-r} \o A_{l+r-n}$ and then induce to 
$A_k\o A_l.$  Denote the corresponding $(A_k \o A_l, A_n \o A_m)$-bimodule 
by $B_r,$ it has a canonical generator that we will call $g_r.$ 
To $r$ there is associated a minimal length representative, $w(r),$ 
of the double 
cosets $\S_k \times \S_l \setminus \S_{k+l} / \S_n \times \S_m.$ 
Namely, $w(r)$ is 
the permutation that preserves the first $r$ elements of the set 
$\{ 1, 2, \dots , n+m\},$ shifts the next $n-r$ elements by 
$k-r$ to the right, shift the following $k-r$ elements by $n-r$ to 
the left and preserves the last $l+r-n$ elements. Sending the 
generator $g_r$ of this bimodule to $Y_{w(r)}\in A_{n+m}$ and summing over 
all admissible $r$ gives us an isomorphism of bimodules 
$\oplusop{r} B_r \cong A_{n+m}.$ Finally, summing 
 over all $(n,m,k,l)$ with $n+m= k+l,$ 
we get a functor isomorphism (\ref{consistency}). 

$\square$ 

\subsection{Coherence relations} 

Bialgebra-categories first appeared in the work of Crane and Frenkel [CF]. 
Crane and Frenkel argued that, while Hopf algebras produce invariants of 
3-manifolds, quantum invariants of 4-manifolds will be governed by 
Hopf categories. In Hopf categories multiplication and comultiplication 
operations become functors, functor isomorphisms take place 
of (co)associativity of (co)multiplication and of the consistency 
relation between multiplication and comultiplication. Crane and Frenkel 
imposed $4$ coherence relations on these functor isomorphisms. 
These relations pop up in our simple example:  

\begin{prop} 
Isomorphisms (\ref{iso-ass}), (\ref{iso-co}) and (\ref{consistency})
satisfy the coherence relations of Crane and Frenkel for bialgebra-categories. 
\end{prop}

The coherence relation for the multiplication can be viewed as 
a cube, depicted in Figure 1. 
 
\drawing{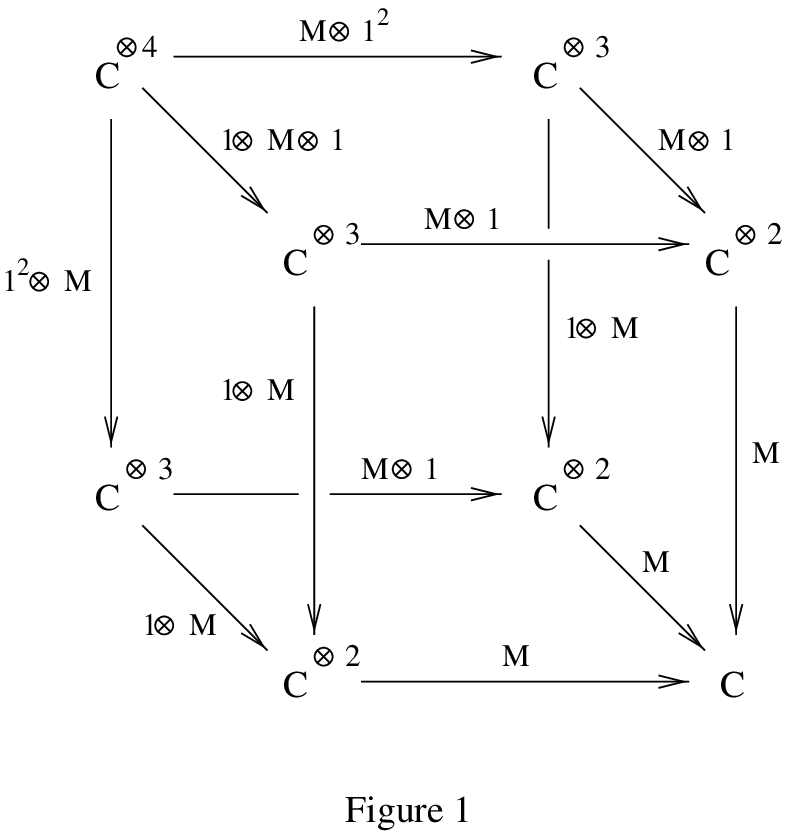} 

In the vertices of the cube we have placed categories, arrows are 
functors and $6$ square facets of the cube are functor isomorphisms. 
For simplicity we write $1$ for the identity functor $\mbox{Id}.$ 
Any path leading from $\C^{\o 4}$ to $\C$ defines a functor, and 
any square facet defines an isomorphism of functors. 
Starting with the functor corresponding to a path, we 
can apply all $6$ isomorphisms and return to the functor we started with. 
The coherence relation requires this natural transformation of 
functors to be the identity. This relation is obvious in our case. 
 Note that the coherence relation 
for the multiplication is just the coherence relation for the tensor 
product functor in the monoidal categories, also known as the 
pentagon associativity (see Mac Lane [M], for instance). 

The coherence relation for the comultiplication is obtained from Figure 1 
by reversing all arrows and  changing all appearances of $M$ into $\Delta,$ 
This coherence holds in our category for obvious reasons too. Moreover,  
if we start from the multiplication coherence relation and  change 
all functors and functor isomorphisms to their right adjoints, we get 
the coherence relation for the comultiplication. Or, if we start 
from the multiplication coherence relation and pass to left adjoints, 
we again get the comultiplication coherence relation (after 
canceling out all appearances of permutations in left adjoints, since
the left adjoint of $M$ is $S_{12}\Delta$). 

There are two coherence relations that contain both multiplication and 
comultiplication. One of them is depicted in Figure 2. To get the other 
one, change Figure 2 in the following way: exchange $M$ with $\Delta$ 
everywhere, reverse the direction of all arrows and invert the 
order of all compositions of functors, i.e., $M^2 \circ S_{23}$ should 
become $S_{23} \circ \Delta^2.$ 

\drawing{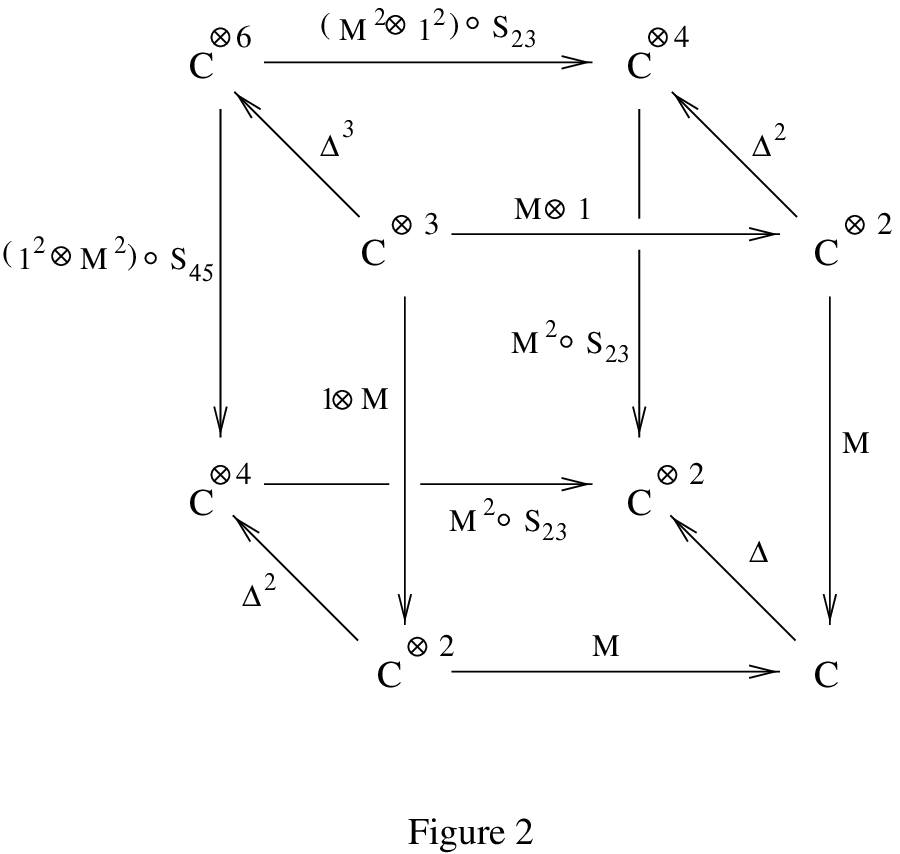} 

In our category $\C$ functors $M$ and $\Delta$ have 
nice adjointness properties, and these two coherence relations 
are equivalent via adjointness.
Figure 2 coherence cube in the category $\C$ follows from a simple 
computation with double cosets $\S_k \times \S_l 
\setminus \S_{n+m+p} / \S_n \times 
\S_m \times \S_p$ for $k+l=n+m+p.$ We omit the details. $\square$

\subsection{Other structures} 

{\bf Commutativity and cocommutativity: } 
The bialgebra $R$ is commutative and cocommutative. The bialgebra-category 
$\C$ is not commutative or cocommutative, in the sense that $MS_{12}$ is 
not isomorphic to $M$ and $S_{12}\Delta$ is not isomorphic to $\Delta.$ 
Instead, we have isomorphisms (\ref{quasi-com-mult}) and 
(\ref{quasi-com-delta}), which say that $MS_{12},$ resp. $S_{12}\Delta$ 
is isomorphic to $M,$ resp. $\Delta,$ twisted by the involution functor 
$\Psi.$ We will refer to these properties of $M$ and $\Delta$ as 
quasi-commutativity and quasi-cocommutativity, respectively. 
What are the coherence relations for quasi-commutativity and 
quasi-cocommutativity? First of all, the usual coherence 
cube for  
the associativity and commutativity constraints in symmetric monoidal 
categories  can be twisted by $\Psi$ into the one, depicted 
in Figure 3 (where $\Psi^2$ denotes $\Psi^{\o 2},$ etc.) 

\drawing{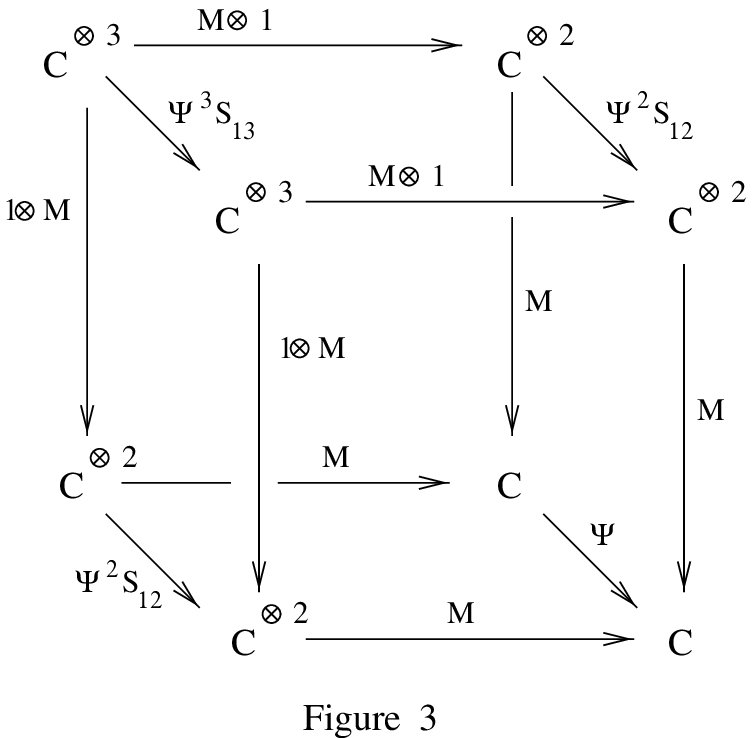} 

To obtain  the coherence relation
between the quasi-cocommutativity and coassociativity isomorphisms,
change the direction of all arrows in Figure 3, substitute $\Delta $ for 
$M$ and invert the order of all compositions. 
Finally, the Figure 4 below shows a coherence cube for the ``mixed'' 
quasi-(co)commutativity.

\drawing{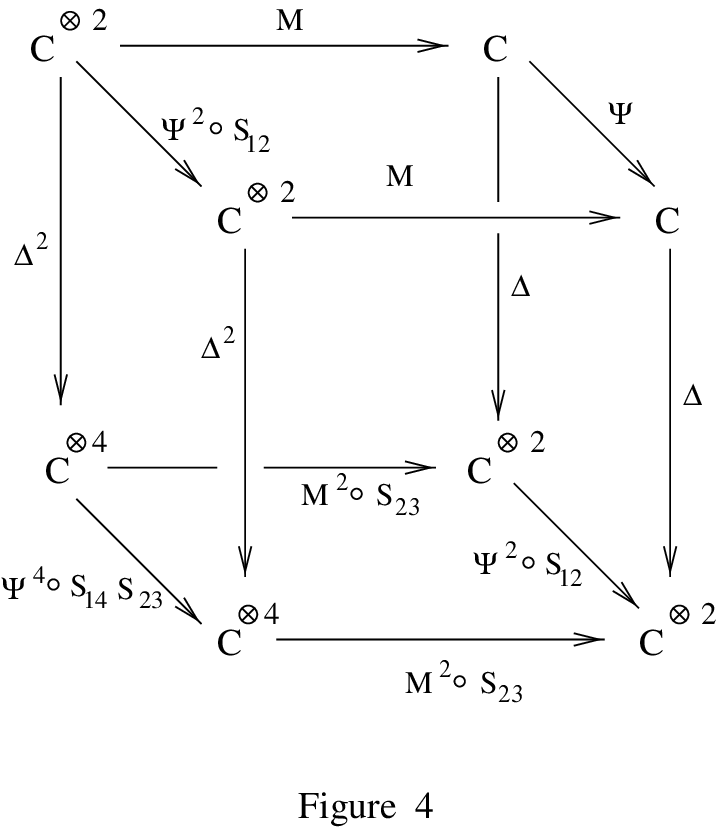} 

\begin{prop} These $3$ coherence relations hold in the category $\C.$ 
\end{prop} 

{\bf Unit and counit: }
Let $\Q\mbox{-vect}$ be the category of finite-dimensional $\Q$-vector 
spaces. The functor $\iota: \Q\mbox{-vect} \to \C$ which takes a vector space 
$V$ to itself, considered as a module over $A_0 =\Q,$ plays the role of 
the unit in the bialgebra-category  $\C.$ 
The functor $\epsilon: \C\to \Q\mbox{-vect}$ 
which takes $\C_n$ to $0$ for $n>0$ and $V\in \C_0$ to $V\in \Q\mbox{-vect}$ 
is the counit functor. 

\vspace{0.2in} 

{\bf Antipode: } 
So far we referred to $\C$ as a bialgebra-category, rather than a Hopf 
category, and did not say a word about the antipode. The antipode 
$s$ in the Hopf algebra $R$ is given by $s(x)= -x.$ Clearly, the antipode 
cannot be lifted to any exact functor in $\C,$ since it does not have 
positive coefficients in the basis of simple modules. This negativity is 
not a serious obstacle, though. We can pass to the bounded derived category 
$D^b(\C)$ of $\C,$ derive the functors $M$ and $\Delta$ and define 
the antipode functor $T: D^b(\C)\to D^b(\C)$ as the composition of 
a shift by $[n]$ (for $\C_n$) 
in the derived category and $\Psi,$
\begin{equation} 
T(N) = \Psi  N [n], \hspace{0.1in } \mathrm{ for }\hspace{0.1in} 
N \in D^b(\C_n).
\end{equation} 
 On the Grothendieck group 
level $\Psi$ does nothing, but it enables   
us to lift the identity $s(ab)= s(b)s(a)$ to the isomorphism of 
functors $ T M \cong  M T^{\o 2} S_{12}.$ But we are in for a bigger 
trouble: there is no functor isomorphism 
\begin{equation} 
\label{trouble} 
M (T \o \mbox{Id})\Delta \cong  \iota\epsilon ,
\end{equation} 
 which any self-respecting 
Hopf category must have. No easy way to save the day by modifying 
the antipode functor is in sight. 
The problem lies with our childish definition of tensor 
powers of $\C,$ as the direct sum of many little blocks. 
One conjectural remedy would be to glue these 
little pieces into a more sophisticated construct,  which should retain 
all the nice bialgebra-category properties of $\C$ and should also
have an antipode functor with the isomorphism (\ref{trouble}) and 
coherence relations for it.

\section{Miscellaneous} 
\label{miss} 

\subsection{Graded bimodules and a categorification of the quantum 
Weyl algebra} \label{graded} 

The algebra $A_n$ is graded, with each $Y_i$ in degree $1,$ and  
the Poincare polynomial of $A_n$ is $[n]!$ where $[n]!= [1]\dots [n]$ and 
$[i]= 1+q+ \dots + q^{i-1}.$ Let $\underline{\C}_n$ be the category of 
finite-dimensional graded left $A_n$-modules and 
$\underline{\C}= \oplusop{n\ge 0}\underline{\C}_n.$ Let $\{i \}$ be the 
functor of shifting the grading up by $i.$ Bimodules $X$ and $D$ over 
$A$ are graded, tensoring with these bimodules over $A$ give us functors, 
denoted $\underline{F}_X,\underline{F}_D,$ in the category $\underline{\C}.$

\begin{prop} There is a functor isomorphism 
\begin{equation} 
\underline{F}_D \underline{F}_X \cong \underline{F}_X\underline{F}_D 
\{ 1\} \oplus \mathrm{Id}
\end{equation} 
\end{prop} 

We define the quantum Weyl algebra as the algebra over $\Z[q,q^{-1}],$ 
generated by $x$ and $\partial,$ with relation $\partial x = qx\partial+ 1.$ 
Let $\underline{R}$ be the module over the quantum Weyl algebra, spanned 
over $\Z[q,q^{-1}]$ 
by $\frac{x^n} {[n]!},$ with the action $x\cdot x^ i = x^{i+1}, 
\partial x^i = [i] x^{i-1}.$ 

The Grothendieck group $K(\underline{\C})$ 
of $\underline{\C}$ is a free $\Z[q,q^{-1}]$-module, 
spanned by the images of simple modules $L_n.$ The $\Z[q,q^{-1}]$-module 
structure comes from the grading, $[N\{ i\} ] = q^i [N],$ for a graded 
module $N.$ Thus, 
\begin{equation} 
[L_n] = \frac{x^n}{[n]!}, \hspace{0.14in} [P_n] = x^n.
\end{equation} 
As a result, $K(\underline{\C})$ can be naturally identified with 
$\underline{R}.$ All other structures described in Section~\ref{Weyl} 
have their graded versions. We skip the details. 

\vspace{0.1in} 

The product $\underline{M}: \underline{\C}^{\o 2}\to \underline{\C}$ 
is again defined as the induction functor, while the coproduct 
$\underline{\Delta},$ considered as a functor from $\underline{\C}_n$ to 
$\oplusop{0\le k\le n}\underline{\C}_k \o \underline{\C}_{n-k},$ 
is the restriction from $A_n$ to $A_{k}\o A_{n-k},$ composed 
with the shift in the grading up by $n-k.$ On the Grothendieck group, 
the coproduct functor acts as the comultiplication
$\Delta(x) = x\o 1 + q\o x.$ 
Functor isomorphisms of Proposition~\ref{funct-iso} hold in the graded 
case as well and all results of Section~\ref{bialgebra} generalize 
easily to the graded case. 

\subsection{Representations of symmetric groups} 

The bialgebra-category $\C$ is reminiscent of 
a similar structure for symmetric groups, discovered by 
Geissinger [G], who observed that  
induction and restriction functors associated to 
inclusions of symmetric groups $\S_n\times \S_m\hookrightarrow \S_{n+m}$ 
induce a bialgebra structure on the direct sum of Grothendieck groups 
of the categories of $\S_n$-modules, over all $n.$ 
Geissinger [G] and 
Zelevinsky [Z] consistently derived many classical results on 
representations of symmetric groups from this Hopf algebra structure. 
Zelevinsky 
 also generalized this construction from symmetric groups to 
wreath products of symmetric groups with 
finite groups and to $GL(n,\mathbb{F}),$ for a finite field $\mathbb{F}.$ 
Although Gessinger and Zelevinsky work mostly with Grothendieck 
groups, their results can be immediately reformulated in terms of categories. 
In particular, induction and restriction define a bialgebra-category 
structure on 
the category $\oplusop{n\ge 0} k[\S_n]\mbox{-mod},$ where 
$k$ is a field and $k[\S_n]\mbox{-mod}$ 
the category of finite-dimensional modules over the group algebra of 
$\S_n.$ 
There are several other interesting examples of bialgebra-categories  
that naturally appear in representation theory. We will discuss 
them elsewhere. 

\subsection{Nil wreath products} 
\label{wreath} 

Let $B$ be an algebra over $\Q.$ By $A_n(B)$ we denote the semidirect 
product of $A_n$ and $B^{\o n},$ with the multiplication 
$ Y_w b_1\o \dots \o b_n = b_{w(1)}\o \dots \o b_{w(n)} Y_w$ for 
$w\in \S_n$ and $b_i\in B.$ 

\begin{prop} If $B$ is a Frobenius algebra then $A_n(B)$ is also Frobenius. 
\end{prop} 

If $B= \Q[z]/\{ z^k =0\},$ denote $A_n(B)$ by $A_n(k).$
 Let $\C_n(k)$ be the category of 
finite dimensional $A_n(k)$-modules and $\C(k)=\oplusop{n\ge 0}\C_n(k).$ 
Inclusions $A_n(k) \hookrightarrow A_{n+1}(k)$ induce induction and 
restriction functors between categories $\C_n(k)$ and $\C_{n+1}(k).$ Denote 
by $F_{X,k},$ resp. $F_{D,k},$ the direct sum of these induction, resp. 
restriction functors over all $n.$ 

\begin{prop} There is a functor isomorphism 
$F_{D,k} F_{X,k} \cong F_{X,k} F_{D,k} \oplus 
(\mathrm{Id}^{\oplus k}).$ 
\end{prop}
Various constructions and results of previous sections, including 
adjointness isomorphisms and bialgebra-category structures, 
can be generalized to algebras $A_n(k).$ 
These algebras are nilpotent counterparts  
of the wreath products of symmetric groups with cyclic groups and 
 of Ariki-Koike cyclotomic Hecke algebras [AK].
For instance, the nilCoxeter and Hecke algebras belong to a two-parameter 
family of algebras with generators $T_1, \dots T_{n-1}$ and relations 
$T_i^2= a T_i + b, T_i T_j = T_j T_i $ for $|i-j|>1 $ and 
$T_i T_{i+1} T_i = T_{i+1} T_i T_{i+1}$ (specializing   $a=b=0, $  resp. 
$a=1-q, b=q,$ gets us the nilCoxeter algebra, resp. the Hecke algebra).   
From this point of view, 
our categorification of the Weyl algebra action on polynomials 
is a toy degeneration of Ariki's magnificent 
realization [A] of irreducible 
highest weight modules over the  affine Lie algebra 
$\widehat{\mf{sl}}_n$ 
as Grothendieck groups of categories of 
modules over cyclotomic Hecke algebras. 
%
%
%

\end{document}